\documentclass[12pt]{amsart}
\usepackage[centertags]{amsmath}
\usepackage{amsfonts,amssymb}
\usepackage{graphicx}
\usepackage{enumerate}

  \newtheorem{theorem}{Theorem}

  \theoremstyle{remark}

  \newcommand{\ZZ}{\mathbb Z}
\begin{document}

\title[A note on Carmichael numbers in residue classes]{A note on Carmichael numbers\\in residue classes}

\date{\today}

\author{Carl Pomerance}
\address{Mathematics Department, Dartmouth College, Hanover, NH 03784}
\email{carl.pomerance@dartmouth.edu}

\begin{abstract}
Improving on some recent results of Matom\"aki and of Wright, we show that
the number of Carmichael numbers to $X$ in a coprime residue class exceeds
$X^{1/(6\log\log\log X)}$ for all sufficiently large $X$ depending on the modulus
of the residue class.
\end{abstract} 
\maketitle
\vskip-30pt
\newenvironment{dedication}
        {\vspace{6ex}\begin{quotation}\begin{center}\begin{em}}
        {\par\end{em}\end{center}\end{quotation}}
\begin{dedication}
{In memory of Ron Graham (1935--2020)\\ and Richard Guy (1916--2020)}
\end{dedication}
\vskip20pt

\section{Introduction}

The ``little theorem" of Fermat asserts that when $p$ is a prime number,
we have $b^p\equiv b\pmod p$ for all integers $b$.  Given two integers
$b,p$ with $p>b>0$, it is computationally easy to check this congruence,
taking $O(\log p)$ arithmetic operations in $\ZZ/p\ZZ$.  So, if the congruence
is checked and we find that $b^p\not\equiv b\pmod p$ we immediately deduce that
$p$ is composite.  Unfortunately there are easily found examples where
$n$ is composite and the Fermat congruence holds for a particular $b$.
For example it always holds when $b=1$.  It holds when $b=2$ and $n=341$,
and another example is $b=3$, $n=91$.  There are even composite numbers
$n$ where $b^n\equiv b\pmod n$ holds for all $b$, the least example being $n=561$.  These are the {\it Carmichael numbers}, named after R. D. Carmichael
who published the first few examples in 1910, see \cite{C}.  (Interestingly,
\v{S}imerka
published the first few examples 25 years earlier, see \cite{S}.)

We now know that there are infinitely many Carmichael numbers, see \cite{AGP},
the number of them at most $X$ exceeding $X^c$ for a fixed $c>0$ and $X$ sufficiently large.

A natural question is if a given residue class contains infinitely many Carmichael numbers.
After work of Matom\"aki \cite{M} and Wright \cite{W1}, we now know
there are infinitely many in a coprime residue class.
 More precisely, we have the following two theorems.
Let 
\[
C_{a,M}(X) =\#\{n\le X: n\hbox{ is a Carmichael number},\,n\equiv a\kern-5pt\pmod{M}\}.
\]

\medskip

{\bf Theorem M} (Matom\"aki). {\it Suppose that $a,M$ are positive coprime
integers and that $a$ is a quadratic residue} mod $M$.  {\it Then
$C_{a,M}\ge X^{1/5}$ for $X$ sufficiently large depending on the choice of $M$.}

\medskip

{\bf Theorem W} (Wright). {\it Suppose that $a,M$ are positive coprime integers.
There are positive numbers $K_M,X_M$ depending on the choice of $M$
such that $C_{a,M}(X)\ge X^{K_M/(\log\log\log X)^2}$ for all $X\ge X_M$.}
\medskip

Thus, Wright was able to remove the quadratic residue condition in Matom\"aki's
theorem but at the cost of lowering the count to an expression that is of the
form $X^{o(1)}$.  The main contribution of this note is to somewhat strengthen
Wright's bound.

\begin{theorem}
\label{thm:main}
Suppose that $a,M$ are positive coprime
integers.   Then $C_{a,M}(X)\ge X^{1/(6\log\log\log X)}$ for all sufficiently large $X$ 
depending on the choice of $M$.
\end{theorem}

That is, we reduce the power of $\log\log\log X$ to the first power and we remove the
dependence on $M$ in the bound, though there still remains the condition
that $X$ must be sufficiently large depending on $M$.  (It's clear though
that such a condition is necessary since if $M>X$ and $a=1$, then
there are no Carmichael numbers $n\le X$ in the residue class $a$~mod~$M$.)

Our proof largely follows Wright's proof of Theorem W, but with a few differences.

Unlike with primes, it is conceivable that a non-coprime residue class contains
infinitely many Carmichael numbers, e.g., there may be infinitely many that
are divisible by 3.  This is unknown, but seems likely.  In fact, it
is conjectured in \cite{BP} that if $\gcd(g,2\varphi(g))=1$, where $g=\gcd(a,M)$,
then there are infinitely many Carmichael numbers $n\equiv a\pmod M$.  
Though we don't know this for any example with $g>1$, 
the old heuristic of Erd\H os \cite{E} suggests that $C_{a,M}(X)\ge X^{1-o(1)}$ as
$X\to\infty$.   

\section{Proof of Theorem \ref{thm:main}}
There is an elementary and easily-proved criterion for Carmichael numbers:
a composite number $n$ is one if and only if it is squarefree and $p-1\mid n-1$
for each prime $p$ dividing $n$.  This is due to Korselt, and perhaps others,
and is over a century old.  In our construction we will have a number $L$
composed of many primes, a number $k$ coprime to $L$ that is not much
larger than $L$, and primes $p$ of the form $dk+1$ where $d\mid L$.
We will show there are many $n\equiv a\pmod M$ that are
 squarefree products of the $p$'s and are $1\pmod{kL}$.  Such $n$, if
 they involve more than a single $p$, will satisfy Korselt's criterion and
 so are therefore Carmichael numbers.

We may assume that $M\ge2$.  Let $\mu=\varphi(4M)$, so that $4\mid\mu$.
Let $y$ be an independent variable; our
other quantities will depend on it.  For a positive integer $n$ let $P(n)$ denote the
largest prime factor of $n$ (with $P(1)=1$), and let $\omega(n)$ denote the
number of distinct prime factors of~$n$.

Let 
\[
{\mathcal Q}_0=\{q \hbox{ prime}:y<q\le y\log^2y,\,q\equiv-1\kern-5pt\pmod{\mu},\,P(q-1)\le y\}.
\]
If $q\le y\log^2y$ and $P(q-1)>y$, then $q$ is of the form $mr+1$, where
$m<\log^2y$ and $r$ is prime.  By Brun's sieve (see \cite[(6.1)]{HR}), the number of such primes $q$ is at most
\[
\sum_{m<\log^2y}\sum_{\substack{r~{\rm prime}\\mr\le y\log^2y\\rm+1~{\rm prime}}}1\ll\sum_{m<\log^2y}\frac{y\log^2y}{\varphi(m)\log^2y}\ll y\log\log y.
\]
Also, the number of primes $q\le y\log^2y$ with $q\equiv-1\pmod{\mu}$
is $\sim\frac1{\varphi(\mu)}y\log y$ as $y\to\infty$ by the prime number
theorem for residue classes.  We conclude that
\begin{equation}
\label{eq:Qcount}
\#{\mathcal Q}_0\sim\frac1{\varphi(\mu)}y\log y
~\hbox{ and }~\prod_{q\in{\mathcal Q}_0}q=\exp\Big(\frac{1+o(1)}{\varphi(\mu)}y\log^2y\Big),\quad y\to\infty.
\end{equation}
We also record that
\begin{equation}
\label{eq:Qrecip}
\sum_{q\in{\mathcal Q}_0}\frac1q=o(1),\quad y\to\infty,
\end{equation}
since this holds for all of the primes in the interval $(y,y\log^2y]$.

Fix $0<B<5/12$; we shall choose a numerical value for $B$ near to
$5/12$ at the end of the argument.  Let
\begin{equation}
\label{eq:x}
x=M^{1/B}\prod_{{q\in \mathcal Q}_0}q^{1/B}.
\end{equation}
It follows from  \cite[(0.3)]{AGP}
that there is an absolute constant $D$ and a set ${\mathcal D}(x)$ of
at most $D$ integers greater
than $\log x$, such that if $n\le x^{B}$, $n$ is not
divisible by any member of ${\mathcal D}(x)$, $b$ is coprime to $n$,
and $z\ge nx^{1-B}$,
then the number of primes
$p\le z$ with $p\equiv b\pmod n$ is $>\frac12\pi(z)/\varphi(n)$.

For each number in ${\mathcal D}(x)$ we choose a prime factor and remove
this prime from ${\mathcal Q}_0$ if it happens to be there.  Let $L$ be the
product of the primes in the remaining set ${\mathcal Q}$, so that $L$
 is not divisible by any member of ${\mathcal D}(x)$, and
$\mathcal Q$ satisfies \eqref{eq:Qcount} and \eqref{eq:Qrecip}.  In particular, 
\begin{equation}\begin{aligned}
\label{eq:Lfacts}
L=&\exp\Big(\frac{1+o(1)}{\varphi(\mu)}y\log^2y\Big),\quad \omega(L)\sim\frac1{\varphi(\mu)}y\log y,\\
&\hbox{and }~\sum_{q\mid L}\frac1q=o(1)
~\hbox{ as }~y\to\infty.
\end{aligned}\end{equation}
In addition, we have $ML\le x^{B}$.

For each $d\mid L$ and each quadratic residue $b\pmod{L/d}$
we consider the primes 
\begin{enumerate}
\item[$\bullet$]
$p\le dx^{1-B}$,
 \item[$\bullet$]
$p\equiv a\pmod M$,
\item[$\bullet$]
 $p\equiv1\pmod d$,
 \item[$\bullet$]
$p\equiv b\pmod{L/d}$.  
\end{enumerate}

Since $M$ is coprime to $L$, the congruences may be glued to a single
congruence modulo $ML$, 
and the number of such primes $p$ is
\[
>\frac{\pi(dx^{1-B})}{2\varphi(ML)}>\frac{dx^{1-B}}{3\varphi(ML)\log x}
\]
for $y$ sufficiently large.

We add these inequalities over the various choices of $b$, the number of which
is $\varphi(L/d)/2^{\omega(L/d)}$, so the number of primes $p$ corresponding
to $d\mid L$ is
\[
>\frac{dx^{1-B}2^{\omega(d)}}{3\cdot2^{\omega(L)}\varphi(Md)\log x}.
\]
We wish to impose an additional restriction on these primes $p$, namely
that $\gcd((p-1)/d,L)=1$.  For a given prime $q\mid L$ the number of 
primes $p$ just counted and for which $q\mid (p-1)/d$ is, via the
Brun--Titchmarsh inequality,
\[
\ll\frac{dx^{1-B}2^{\omega(d)}}{2^{\omega(L)}q\varphi(Md)\log(x/(qML))}
\ll\frac{dx^{1-B}2^{\omega(d)}}{2^{\omega(L)}q\varphi(Md)\log x}.
\]
Summing this over all $q\mid L$ and using that $\sum_{q\mid L}1/q=o(1)$,
these primes $p$ are seen to be negligible.  It follows that for $y$ sufficiently
large, there are 
\[
>\frac{dx^{1-B}2^{\omega(d)}}{2^{\omega(L)+2}\varphi(Md)\log x}
>\frac{x^{1-B}2^{\omega(d)}}{2^{\omega(L)+2}\varphi(M)\log x}
\]
primes $p\le dx^{1-B}$ with $p\equiv 1\pmod d$,
$\gcd((p-1)/d,L)=1$, $p\equiv a\pmod M$, and
$p$ is a quadratic residue (mod~$L$) (noting that 1 (mod~$d$) is a quadratic
residue (mod $d$)).

For each pair $p,d$ as above, we map it to $(p-1)/d$ which is an integer
$\le x^{1-B}$ coprime to $L$.  The number of pairs $p,d$
is
\[
>\frac{x^{1-B}}{2^{\omega(L)+2}\varphi(M)\log x}\sum_{d\mid L}2^{\omega(d)}
=\frac{x^{1-B}3^{\omega(L)}}{2^{\omega(L)+2}\varphi(M)\log x}.
\]
We conclude that there is a number $k\le x^{1-B}$ coprime to $L$ which has 
$>(3/2)^{\omega(L)}/(4\varphi(M)\log x)$
representations as $(p-1)/d$.  Let $\mathcal P$ be the set of primes $p=dk+1$
that arise in this way.  Then
\begin{equation}
\label{eq:Pcount}
\#{\mathcal P}>\frac{(3/2)^{\omega(L)}}{4\varphi(M)\log x}.
\end{equation}

For a finite abelian group $G$, let $n(G)$ denote Davenport's constant,
 the least number
such that in any sequence of group elements of length $n(G)$ there is
a non-empty subsequence with product the group identity.  It is easy to
see that $n(G)\ge\lambda(G)$ (the universal exponent for $G$), and in
general it is not much larger: $n(G)\le\lambda(G)(1+\log(\#G))$.  This result
is essentially due to van Emde Boas--Kruyswijk and Meshulam, see \cite{AGP}.

Let $G$ be the subgroup  of $(\ZZ/kML\ZZ)^*$ of residues 
$\equiv1\pmod k$.  We have $\#G\le ML$.  Also, $\lambda(G)
\le M\lambda(L)$.  (Note that, as usual, we denote $\lambda((\ZZ/L\ZZ)^*)$
by $\lambda(L)$.  It is the lcm of $q-1$ for primes $q\mid L$, using that $L$ is
squarefree.)  Each prime dividing $\lambda(L)$ is at most $y$
and each prime power dividing $\lambda(L)$ is at most $y\log^2y$, so that
\[
\lambda(L)\le(y\log^2y)^{\pi(y)}.
\]
Thus, for large $y$, using \eqref{eq:Lfacts},
\begin{equation}
\label{eq:Dav}
n(G)\le M(y\log^2y)^{\pi(y)}\log(ML)\le e^{2y}.
\end{equation}

For a sequence $A$ of elements in a finite abelian group $G$, let $A^*$
denote the set of nonempty subsequence products of $A$.  In 
Baker--Schmidt \cite[Proposition 1]{BS} it is shown that there is a number $s(G)$ such that if $\#A\ge s(G)$,
 then $G$ has
a nontrivial subgroup $H$ such that $(A\cap H)^*=H$.  Further,
\[
s(G)\le 5\lambda(G)^2\Omega(\#G) \log(3\lambda(G)\Omega(\#G)),
\]
where $\Omega(m)$ is the number of prime factors of $m$ counted with multiplicity.
Thus, with $G$ the group considered above, we have
\[
s(G)\le e^{3y}
\]
for $y$ sufficiently large.  

It is this theorem that Matom\"aki and Wright use in their papers on Carmichael
numbers.  The role of the sequence $A$ is played by $\mathcal P$, the set
of primes constructed above of the form $dk+1$ where $d\mid L$.  So, if
$\#{\mathcal P}>s(G)$ we are guaranteed that every member of a nontrivial
subgroup $H$ of $G$ is represented by a subset product of ${\mathcal P}\cap H$.

We don't know precisely what this subgroup $H$ is, but we do know that it is nontrivial
and that it is generated by members of $\mathcal P$.  Well, suppose $p_0$ is in
${\mathcal P}\cap H$.  Then $p_0^m\in H$ for every integer $m$.
Note that by construction, $\gcd(\lambda(L)/2,\varphi(M))=1$, so there
is an integer $m\equiv 1\pmod{\varphi(M)}$ and $m\equiv0\pmod{\lambda(L)/2}$.
Further, since $p_0$ is a quadratic residue (mod $L$), it follows that
$p_0^{\lambda(L)/2}\equiv1\pmod{L}$.  Thus, $p_0^m\equiv1\pmod L$
and $p_0^m\equiv a\pmod M$ (since $m\equiv1\pmod{\varphi(M)}$.

Thus, there is a subsequence product $n$ of $\mathcal P$ that is
1 (mod $kL$) and $a$ (mod $M$).  (Note that every member of $G$ is 1 (mod $k$.)
Further, $n$ is squarefree and for each prime factor $p$ of $n$ we have $p-1\mid kL$.
Since $n\equiv1\pmod{kL}$ we have $p-1\mid n-1$.  Thus, $n\equiv a\pmod M$
is either a prime or a Carmichael number.

We actually have many subsequence products $n$ of $\mathcal P$ that satisfy
these conditions, and $\mathcal P$ has at most one element that is 1 (mod $L$),
so we do not need to worry about the case that $n$ is prime.  We let
$t=\lceil e^{3y}\rceil$, so that $t\ge s(G)$.  As shown in \cite{M},
\cite{W1}, the Baker--Schmidt result implies that ${\mathcal P}$
has at least
\[
N:=\binom{\#{\mathcal P}-n(G)}{t-n(G)}\Big/\binom{\#{\mathcal P}-n(G)}{n(G)}
\]
subsequence products $n$ of length at most $t$ which are Carmichael numbers
in the residue class $a\pmod M$.  Thus,
\begin{align*}
N&>\left(\frac{\#{\mathcal P}-n(G)}{t-n(G)}\right)^{t-n(G)}(\#{\mathcal P})^{-n(G)}\\
&>\left(\frac{\#{\mathcal P}}{t}\right)^{t-n(G)}(\#{\mathcal P})^{-n(G)}
>(\#{\mathcal P})^{t-2n(G)}t^{-t}.
\end{align*}

Let $X=x^t$.  Since each $p\in\mathcal P$ has $p\le x$, it follows that
all of the Carmichael numbers constructed above are at most $X$.
Using \eqref{eq:Qcount}, \eqref{eq:x}, and \eqref{eq:Dav}, we have
\[
X=\exp\Big(\frac{1/B+o(1)}{\varphi(\mu)}ty\log^2y\Big),
\]
and using \eqref{eq:Pcount} and \eqref{eq:Lfacts} gives
\begin{align*}
N&\ge\exp\Big(\frac{\log(3/2)+o(1)}{\varphi(\mu)}ty\log y-t\log t\Big)\\
&=\exp\Big(\frac{\log(3/2)+o(1)}{\varphi(\mu)}ty\log y\Big).
\end{align*}
Thus, $N\ge X^{(B\log(3/2)+o(1))/\log y}$.  Now, 
\[
\log X\sim \frac1{B\varphi(\mu)}ty\log^2y,
\]G
so that using $t=\lceil e^{3y}\rceil$,
\[
\log\log X=3y+O(\log y),\quad \log\log\log X = \log y+O(1).
\]
We thus have $N\ge X^{(B\log(3/2)+o(1))/\log\log\log X}$.  The number
$B<5/12$ can be chosen arbitrarily close to $5/12$ and since $(5/12)\log(3/2)>1/6$,
the theorem is proved.

%\noindent{\bf Acknowledgments}.

\end{document}